\documentclass[11pt]{article}
\usepackage{latexsym}
\usepackage[dvips]{epsfig}

\newcommand{\Pc}{{\cal P}}

\newcommand{\Q}{{\cal V}}
\newcommand{\qed}{\mbox{$\Diamond$}\vspace{\baselineskip}}

\newtheorem{theorem}{Theorem}[section]
\newtheorem{proposition}[theorem]{Proposition}
\newtheorem{lemma}[theorem]{Lemma}
\newtheorem{definition}[theorem]{Definition}
\newtheorem{corollary}[theorem]{Corollary}
\newtheorem{example}[theorem]{Example}

\newenvironment{proof}{\noindent {\bf Proof:}}{{\qed}}

\newcommand{\vanish}[1]{} 
\begin{document}
\title{A combinatorial proof of the log-concavity of the numbers of
permutations with $k$ runs}

\author{Mikl\'os B\'ona
\thanks{School of Mathematics,
Institute for Advanced
Study, Princeton, NJ 08540. Supported by Trustee Ladislaus von Hoffmann,
 the Arcana Foundation. Email: bona@math.ias.edu.}
\and Richard Ehrenborg \thanks{
School of Mathematics, Institute for Advanced Study, Princeton, NJ 08540. 
Supported by National Science Foundation,
DMS 97-29992, and NEC Research Institute, Inc.
Email: jrge@math.ias.edu.}}

 \date{}

\maketitle

\begin{abstract}
We combinatorially prove that the number $R(n,k)$ of permutations of
 length $n$ having $k$
runs is a log-concave sequence in $k$, for all $n$. We also give a new 
combinatorial proof for the log-concavity of the Eulerian numbers.
\end{abstract}
\section{Introduction}
\begin{figure}[ht]
 \begin{center}
  \epsfig{file=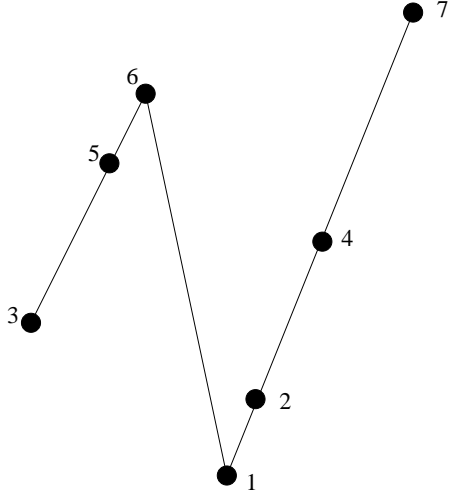}
  \caption{The permutation 3561247 has three runs }
  \label{threeruns}
 \end{center}
\end{figure}

Let $p = p_1 p_2 \cdots p_n$ be a permutation of the set
$\{1,2,\ldots,n\}$ written in the one-line notation. We say that $p$
changes direction at position $i$, if either $p_{i-1}<p_i>p_{i+1}$, or
$p_{i-1}>p_i<p_{i+1}$, in other words, when $p_i$ is either a {\em
peak} or a {\em valley}. We say that $p$ has $k$ runs if there are
$k-1$ indices $i$ so that $p$ changes direction at these positions.
So for example, $p=3561247$ has 3 runs as $p$ changes direction when
$i=3$ and when $i=4$. A geometric way to represent a permutation and
its runs by a diagram is shown on Figure \ref{threeruns}.  The runs
are the line segments (or edges) between two consecutive entries where
$p$ changes direction. So a permutation has $k$ runs if it can be
represented by $k$ line segments so that the segments go ``up'' and
``down'' exactly when the entries of the permutation do. The theory of
runs has been studied in~\cite[Section 5.1.3]{knuth} in connection
with sorting and searching.

In this paper, we are going to study the numbers $R(n,k)$ of permutations
 of length $n$ or, in what follows, $n$-permutations with $k$ runs. We will
show that for any fixed $n$, the sequence $R(n,k)$, $k=0,1,\ldots, n-1$ is
log-concave, that is, $R(n,k-1) \cdot R(n,k+1)\leq R(n,k)^2$.
 In particular, this implies \cite{brenti,stanley} 
that this same sequence is unimodal, that is, there exists an $m$
so that $R(n,1)\leq R(n,2)\leq \cdots \leq R(n,m)\geq R(n,m+1)\geq \cdots
\geq R(n,n-1)$. We will also show that
roughly half of the roots of the generating function $R_n(x)=
\sum_{k=1}^{n-1}R(n,k)x^{k}$ are equal to $-1$, and give a combinatorial
interpretation for the term which remains after one divides $R_n(x)$ by
all the $(x+1)$ factors. While doing that, we will also give a new proof of
the well-known fact \cite{vesco,wilf} that the Eulerian numbers are 
log-concave.

\section{The Factorization of $R_n(x)$} 

Let $p=p_1p_2\cdots p_n$ be a permutation. We say that $i$
is a {\em descent} of $p$ if $p_i>p_{i+1}$, while we say that 
 $i$ is an {\em ascent} of $p$ if $p_i<p_{i+1}$. 

In our study of $n$-permutations with a given number of runs,
 we can clearly assume that 1 is an ascent of $p$.
 Indeed, taking the permutation $q=q_1q_2\cdots q_n$, where $q_i=n+1-p_i$,
 we get the {\em complement} of $p$, which has the same number of runs
as $p$. This implies in particular that for any given $i$, there are as
many $n$-permutations with $k$ runs in which $p_i<p_{i+1}$ as there are
such permutations in which $p_i>p_{i+1}$. Let $Q_n(x)$ be any generating
function enumerating $n$-permutations according to some statistics. We
say that $ Q_n(x)$ is {\em invariant to $(i,i+1)$ } if the set of
$n$-permutations with  $p_i>p_{i+1}$ contributes  $Q_n(x)/2$ to $Q$.
Certainly, in this case the set of $n$-permutations with  $p_i<p_{i+1}$ 
contributes accounts $Q_n(x)/2$ to $Q$ as well.

Let  $R_n(x)=
\sum_{k=1}^{n-1}R(n,k)x^{k}$ be the ordinary generating function of
$n$-permutations with $k$ runs, where $1 \leq k \leq n-1$.
So we have $R_2(x)=2x$,
$R_3(x)=2x+4x^2$, and $R_4(x)=10x^3+12x^2+2x$. One sees that all coefficients
 of $R_n(x)$ are even, which is explained by the symmetry described above.

The following proposition is our initial step of factoring $R_n(x)$. It will
lead us to the definition of an important version of this polynomial.
\begin{proposition}\label{easy} For all $n\geq 4$, the polynomial $R_n(x)$
 is divisible by $(x+1)$. \end{proposition} 
\begin{proof}
It is straightforward to verify 
(by considering all possible patterns for the last four
entries of~$p$) that the involution $I_1$
 interchanging $p_{n-1}$ and $p_n$ increases the number
of runs by 1 in half of all permutations, and, being an involution, decreases
the number of runs by 1 in the other half of the permutations. In particular,
there are as many permutations with an odd number of runs as there are with an
even number of runs, so $(x+1)$ is indeed a divisor of $R_n(x)$.  
\end{proof}

\begin{example} {\em If $n=4$, then there is 1 permutation with 1 run, 6
permutations with 2 runs and 5 permutations with 3 runs. So $R_4(x)=
2(5x^3+6x^2+x)=2(x+1)(5x^2+x)$. }\end{example}
We want to extend the result of Proposition~\ref{easy} by proving that
$R_n(x)$ has  $(x+1)$ as a factor with a large multiplicity, and also, we want
to find a combinatorial interpretation for the polynomial obtained after
dividing $R_n(x)$ by the highest possible power of $(x+1)$. For that purpose,
we introduce the following definition.
\begin{definition}
For $j\leq m=\lfloor (n-2)/2 \rfloor $,
we say that  $p$ is  a {\em $j$-half-ascending permutation} if,
 for all positive
integers $i\leq j$, we have  $p_{n+1-2i} < p_{n+2-2i}$. If $j=m$, then
we will simply say that $p$ is a {\em half-ascending  permutation}.
\end{definition} 
 So $p$ is a 1-half-ascending permutation
permutation if $p_{n-1}<p_n$. In a $j$-half-ascending  permutation, we have $j$
relations, and they involve the rightmost  $j$ disjoint pairs of entries.
The term half-ascending refers to the fact that at least half of the
involved positions are ascents. There are $n!\cdot 2^{-j}$
 $j$-half-ascending permutations
 
Now we define a modified version of the polynomials $R_{n}(x)$
for $j$--half-ascending permutations.
\begin{definition}
Let $p$ be a $j$-half-ascending permutation. Let $r_j(p)$ be the number
 of runs
of the substring $p_1,p_2,\ldots, p_{n-2j}$, and let $s_j(p)$ be the number
of descents of the substring $p_{n-2j},p_{n+1-2j}, \ldots, p_n$. 
Denote $t_j(p)=r_j(p)+s_j(p)$, and define
\[R_{n,j}(x)=\sum_{p\in S_n}x^{t_j(p)}. \]
In particular, we will denote
$R_{n,m}(x)$ by $T_{n}(x)$, that is, $T_{n}(x)$ is the generating
function for half-ascending permutations.
\end{definition} 
So in other words, we count the runs in the non-half-ascending part
 and we count the descents in the half-ascending part (and on that part, as 
it will be discussed, the number of descents determines that of runs.)

\begin{corollary} \label{factor} For all $n\geq 4$ we have
\[\frac{R_n(x)}{x+1}=R_{n,1}(x).\]
Moreover, $R_{n,1}(x)$ is invariant to $(i,i+1)$ for all $i\leq n-3$.
\end{corollary}
\begin{proof}  
Recall from proof of Proposition~\ref{easy} that involution $I_1$ makes
pairs of permutations, and each pair contains two elements whose numbers
of runs differ by 1. Note that half of these pairs consist of two
elements with $p_{n-3}<p_{n-2}$ and the other half consist of two elements with
$p_{n-3}>p_{n-2}$. As $R_n(x)$ is invariant to $(n-3,n-4)$,  it suffices to 
consider the first case.
Dividing $R_n(x)$ by $(x+1)$ we obtain the run-generating
function for the set of permutations which contains one element of each of 
these pairs, namely, the one having the smaller number of runs. 
Observe that for these permutations, the number of runs is equal to the value
 of $t_1(p)$ for the permutation in that pair in which $p_{n-1}<p_n$
 (by checking both possibilities $p_{n-2}<p_{n-1}$ and $p_{n-2}>p_{n-1}$),
so $R_n(x)/(x+1)=R_{n,1}(x)$. Note that our argument also proves that
 those permutations with $p_{i}<p_{i+1}$
 contribute exactly $R_{n,1}(x)/2$ to  $R_{n,1}(x)$, as they represent
half of $R_n(x)$ divided by $(x+1)$, so our second claim is proved, too.
\end{proof}

We point out that it is not true in general that in each pair made by $I_1$,
the permutation having the smaller number of runs is the one with 
$p_{n-1}<p_n$. What is true is that we can {\em suppose} that $p_{n-1}<p_n$
if we count permutations by the defined parameter $t_1(p)$ instead of the
 number of runs. This latter could be viewed as the $t_0(p)$ parameter. 
\begin{example} {\em If $n=4$, then we have 6 permutations in which $p_3<p_4$
and  $p_1<p_2$: 1234, 1324, 1423, 2314, 2413, 3412. We have $t_1(1234)=1$ and
$t_1(p)=2$ for all the other five permutations,  showing
that indeed, $R_{4,1}(x)=2(5x^2+x)$. }\end{example}
For
$1\leq j\leq m$, let $I_j$ be the involution  interchanging
$p_{n+1-2j}$ and $p_{n+2-2j}$.  
Then the following strong result generalizes Proposition~\ref{easy}.  
\begin{lemma} \label{difficult}
For all $n\geq 4$ and $1\leq j \leq m$,  we have
\[\frac{R_n(x)}{(x+1)^j}=R_{n,j}(x),\]
where $m=\lfloor (n-2)/2 \rfloor.$ Moreover, $R_{n,j}(x)$ is invariant to
$(i,i+1)$ for $i\leq n-2j-1$.
\end{lemma}
\begin{proof}
By induction on $j$. For $j=1$, the statement is true by Proposition
\ref{easy}. Now suppose we know the statement for $j-1$.
 
To prove that  $R_{n,j-1}(x)/R_{n,j}(x)=x+1$, we need to group all 
$(j-1)$-half-ascending permutations in pairs,
so that the $t_{j-1}$ values of the
two elements of any given pairs differ by one, and show that
the set of permutations consisting of the elements of each pair having
the smaller $t_{j-1}$ value yields the generating function $R_{n,j}(x)$.

However, $I_j$ just does that, as can be checked by verifying both
 possibilities   $p_{n-2j}<p_{n+1-2j}$ and $p_{n-2j}>p_{n+1-2j}$. These
are the only cases to consider as we can assume by our induction hypothesis
that $p_{n-1-2j}<p_{n-2j}$. Moreover, permutations with $p_{i}<p_{i+1}$
contribute exactly $R_{n,j}(x)/2$ to $R_{n,j}(x)$ if $i\leq n-2j-1$ 
as they represent half of $R_{n,j-1}(x)$ divided by $(x+1)$.
\end{proof}

Note that we have just repeated the proof of Proposition~\ref{easy} with 
general $j$, instead of $j=1$.
\begin{corollary}
We have \[\frac{R_n(x)}{(x+1)^m}=T_{n}(x).\]
\end{corollary}
So we have proved that $m=\lfloor (n-2)/2 \rfloor $ 
of the roots of $R_n(x)$ are equal to 
$-1$, and certainly, one other root is equal to 0 as all permutations have at
least one run. It is possible to prove analytically \cite{wilf} that the
other half of the roots of $R_n(x)$, that is, the roots of $T_{n}(x)$,
 are all, real, negative, and distinct.
That implies~\cite{stanley} that the coefficients of $R_n(x)$ and $T_{n}(x)$
are log-concave.

However, in the next section we will {\em combinatorially} prove that
the coefficients of $T_{n}(x)$ form a log-concave sequence. Let
$U(n,k)$ be the coefficient of $x^k$ in $T_{n}(x)$.  Let ${\cal
U}(n,k)$ be the set of half-ascending permutations with $k$ descents,
so $|{\cal U}(n,k)| = U(n,k)$.

Now suppose for shortness that $n$ is even and assume that $p$ is a
 half-ascending permutation, that is, $p_{2i-1}<p_{2i}$
for all $i$, $1\leq i \leq n/2$. The following proposition
summarizes the different ways we can describe the same parameter of $p$.
\begin{proposition} \label{trivi} Let $p$ be a half-ascending permutation.
 Then $p$ has $2k+1$ runs if and only if $p$ has $k$ descents, or, in other
words, when $t(p)=k+1$. 
\end{proposition}

If $n$ is odd, then the rest of our argument is a little more tedious, though
conceptionally not more difficult. We do not want to break the course of our
proof here, so we will go on with the assumption that $n$ is even, then, in the
second part of the
proof of Theorem~\ref{main}, we will indicate what modifications are necessary
to include the case of odd $n$.

So in order to prove that the sequence $R(n,k)$ is log-concave in $k$, we
need to prove that the  sequence $U(n,k)$ enumerating  half-ascending
$n$-permutations with $k$ descents is log-concave. That would be sufficient as
the convolution of two log-concave sequences is log-concave \cite{stanley}.  

\section{A lattice path interpretation}

Following \cite{vesco}, we will set up a bijection from the set ${\cal A}
(n,k)$ of $n$-permutations with $k$ descents onto that of labeled 
northeastern lattice paths with $n$ edges, exactly $k$ of which are vertical.
However, our lattice paths will be different from those in \cite{vesco};
in particular, they will preserve the information if the position $i$ is
an ascent or descent.

Let $\Pc(n)$
be the set of
labeled northeastern lattice paths with the $n$ edges
$a_1, a_2, \ldots ,a_n$
and the corresponding positive integers as labels
$e_1, e_2, \ldots ,e_n$
so that the following hold:
\begin{enumerate} 
\item[(1)]
the edge $a_1$ is horizontal and $e_1=1$,

\item[(2)]
if the edges $a_i$ and $a_{i+1}$ are both vertical, or both horizontal,
then $e_i\geq e_{i+1}$,

\item[(3)]
if  $a_i$ and $a_{i+1}$ are perpendicular to each other,
then $e_i+ e_{i+1}\leq i+1.$
\end{enumerate}
We will not distinguish between paths which 
can be obtained from each other by translations.
Let $\Pc(n,k)$ be the set of all such labeled lattice paths
which has $k$ vertical edges, and let $P(n,k)=|\Pc (n,k)|$.

\begin{proposition}
\label{properties} 
The following two properties of paths in $\Pc(n)$ are
immediate from the definitions.
\begin{itemize}
\item For all $i\geq 2$, we have $e_i \leq i-1$.
\item Fix the label $e_i$. Then if $e_{i+1}$ can take value $v$, then it 
can take all nonnegative integer values $w\leq v$.
\end{itemize}
\end{proposition}
Also note that all restrictions on $e_{i+1}$ are given by $e_i$, independently
of preceding $e_j$, $j<i$. The following bijection is the main result in this
section.
\begin{theorem}
\label{bij}
The following description defines
a bijection from ${\cal A}(n)$ onto $\Pc(n)$.
Let $p$ be a permutation on $n$ elements.
To obtain the edge $a_{i}$ and the label $e_{i}$
for $2 \leq i \leq n$,
restrict the permutation~$p$ to the $i$ first entries and
relabel the entries to obtain
the permutation $q = q_1 \cdots q_{i}$.
\begin{itemize}
\item
If the position $i-1$ is a descent of the permutation $p$
(equivalently, of the permutation $q$),
let the edge $a_{i}$ be vertical
and the label $e_{i}$ be equal to $q_i$.

\item
If the position $i-1$ is an ascent of the permutation $p$,
let the edge $a_{i}$ be horizontal
and the label~$e_{i}$ be $i + 1 - q_i$.
\end{itemize}
Moreover, this bijection restricts naturally to a bijection
between ${\cal A}(n,k)$ and $\Pc(n,k)$ for $0 \leq k \leq n-1$. 
\end{theorem}
\begin{proof}
It is straightforward to see that the map described is injective on
the set of labeled lattice path, not necessarily satisfying
conditions (2) and (3).
Assume that $i$ and $i+1$ are both descents of the permutation $p$.
Let $q$, respective $r$, be the permutation when restricting to
the $i$, respective $i+1$, first elements. Observe that
$q_i$ is either $r_i$ or $r_i - 1$. Since $r_i > r_{i+1}$
we have $q_i \geq r_{i+1}$ and condition~(2) is satisfied in this case.
By similar reasoning the three remaining cases are shown,
hence the map is into the set $\Pc(n)$.

To see that this is a bijection, we show that we can recover the
permutation $p$ from its image.
It is sufficient to show that we can recover $p_n$, and then
use induction on $n$ for the rest of $p$. To recover $p_n$ from its image,
simply recall that $p_n$ is equal to the label $l$ of the last edge if
that edge is vertical, and to $n+1-l$ if that edge is horizontal.
\end{proof}

\vanish{
We prove the statement by induction on $n$. For $n=1$ the statement is 
trivially true. Suppose we know the statement for $n-1$,
and let $p \in {\cal A}(n,k)$.
Take away the last entry of $p$, this way we obtain a permutation $p'$ which
is either in  ${\cal A}(n-1,k-1)$ or in ${\cal A}(n-1,k)$. (If the entry we 
took away from $p$ was $x$, then we get $p'$ by subtracting 1 from all
 remaining entries larger than $x$, and leaving other entries unchanged).
 In either case, by our induction
hypothesis, $p'$ is mapped into a lattice path either in  $\Pc(n-1,k-1)$
 or
in $\Pc(n-1,k)$. To get the image of $p$ in ${\cal A}(n,k)$, we add 
a vertical edge to
the end of the image of $p'$ in the first case, and a horizontal edge in the
second case. If $n-1$ was a descent of $p$, that is, in the first case, we
label this new vertical edge by $p_n$. If $n-1$ was an ascent of~$p$, that is,
in the second case, we label the new horizontal edge by $n+1-p_n$. One then 
checks easily that the three conditions on the labels are indeed satisfied.

To see that this is a bijection, we show that we can recover $p$ from its 
image. Again, it is sufficient to show that we can recover $p_n$, and then
use induction on $n$ for the rest of $p$. To recover $p_n$ from its image,
simply recall that $p_n$ is equal to the label $l$ of the last edge if
that edge is vertical, and to $n+1-l$ if that edge is horizontal.
}

\begin{figure}[ht]
 \begin{center}
  \epsfig{file=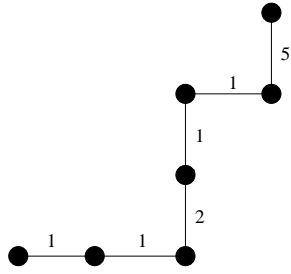}
  \caption{The image of the permutation $243165$}
  \label{onepath}
 \end{center}
\end{figure}

The lattice path corresponding to the permutation $243165$ is shown on 
Figure~\ref{onepath}.

The difference between our bijection and that of \cite{vesco} is that in ours,
the direction of $a_i$ tells us whether $p_{i-1}$ is a descent in $p$.
This is why we can use this bijection to gain information the class of 
half-ascending permutations.
\begin{corollary}
The bijection in Theorem~\ref{bij}
restricts to a bijection from
${\cal U}(n,k)$ to
lattice paths in $\Pc(n,k)$
where $a_i$ is horizontal for all even indices $i$.
\end{corollary}

\section{The log-concavity of $U(n,k)$}
In this section we are going to give a new proof for the fact that the
numbers $A(n,k)=|{\cal A}(n,k)|$ are unimodal in $k$, for any fixed $n$. 
This fact is already
known and has elegant proofs \cite{vesco}. However, our proof will also 
indicate the unimodality of the $U(n,k)$.
\begin{theorem} \label{logcon}
For all positive integers $n$ and all positive integers $k\leq n$ we
have \[A(n,k-1)\cdot A(n,k+1)\leq A(n,k)^2 \] 
and also \[U(n,k-1)\cdot U(n,k+1)\leq U(n,k)^2.\]
\end{theorem}
\begin{proof} To prove the theorem combinatorially, we construct a {\em 
quasi-injection}  
\[\Phi: \Pc(n,k-1) \times \Pc(n,k+1) \longrightarrow
       \Pc(n,k) \times \Pc(n,k).\]
By quasi-injection we mean that there will be some elements of 
$\Pc(n,k-1) \times \Pc(n,k+1)$ for which $\Phi$ will not be defined, 
but the number of these elements will be less than that of elements 
in $\Pc(n,k) \times \Pc(n,k)$ which are not in the image of $\Phi$. 

In particular, the restriction of $\Phi$ onto $\Q(n,k-1) \times \Q(n,k+1)$
will map into $\Q(n,k) \times \Q(n,k)$, where $\Q(n,k)$ is the subset of
$\Pc(n,k)$ consisting of lattice paths in which $a_i$ is horizontal for all
even~$i$. 
Let $(P,Q)\in \Pc(n,k-1) \times \Pc(n,k+1)$. Place the initial points of
 $P$ and
$Q$ at $(0,0)$ and $(1,-1)$, respectively. Then the endpoints of $P$ and
$Q$ are $(n-k+1,k-1)$ and $(n-k,k)$, respectively, so $P$ and $Q$ intersect.
Let $X$ be their {\em first} intersection point (we order intersection points
from southwest to northeast), and decompose $P=P_1\cup P_2$ and 
$Q=Q_1\cup Q_2$, where $P_1$ is a path from $(0,0)$ to $X$, 
 $P_2$ is a path from $X$ to $(n-k,k)$, $P_1$ is a path from  $(1,-1)$ to $X$,
and $Q_2$ is a path from $X$ to $(n-k+1,k-1)$. Let $P'=P_1 \cup Q_2$ and
let  $Q'=Q_1 \cup P_2$. If $P'$ and $Q'$ are valid paths, that is, if their
labeling fulfills conditions (1)--(3), then we set $\Phi(P,Q)=(P',Q')$. 
See Figure \ref{twopaths} for this construction.

\begin{figure}[ht]
 \begin{center}
  \epsfig{file=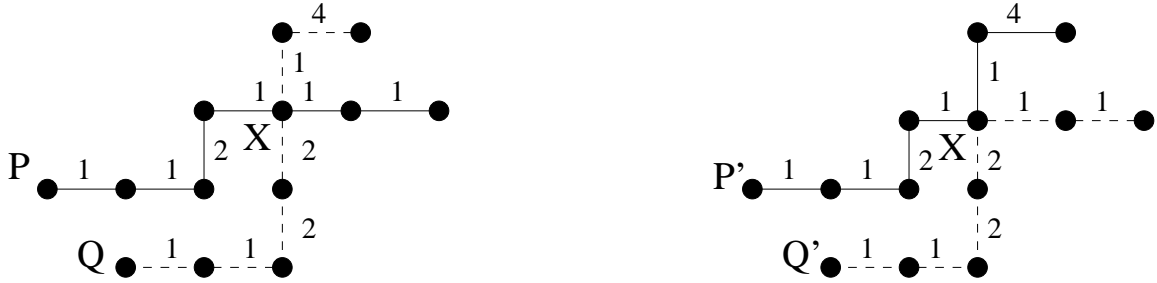}
  \caption{ Constructing the new pair of paths}
  \label{twopaths}
 \end{center}
\end{figure}

It is clear that $\Phi(P,Q)=(P',Q')\in  \Pc(n,k) \times \Pc(n,k)$,
 (in particular,
 $(P',Q')$ belongs to the subset of $\Pc(n,k) \times \Pc(n,k)$ consisting of 
{\em intersecting} pairs of paths), and that $\Phi$ is one-to-one. 
What remains  to show is that the number of pairs $(P,Q)\in \Pc(n,k-1) \times
 \Pc(n,k+1)$
for which $\Phi$ cannot be defined this way is less than the number of pairs
$(P',Q')\in  \Pc(n,k) \times \Pc(n,k)$ which are not obtained as
 images of $\Phi$.
In fact, we will show that this is true even if we restrict ourselves to
pairs $(P',Q')\in  \Pc(n,k) \times \Pc(n,k)$ which do intersect.

Let $a,b,c,d$ be the labels of the 
four edges adjacent to $X$ as shown in Figure \ref{cross}, the edges $AX$ and
$XB$ originally belonging to $P$ and the edges $CX$ and $XD$ originally
 belonging to $Q$.
 (It is possible that these four edges are not all distinct; $A$ and $C$ are
always distinct as $X$ is the first intersection point, but it could be,
 that $B=D$ and so $BX=DX$; this singular case can be treated very
 similarly to the generic case we describe below and 
hence omitted). Then 
 a configuration shown on Figure \ref{cross}
 can be part of a pair $(P,Q)$ in the
domain of $\Phi$  exactly when $a\geq b$ and  $c\geq d$.
On the other hand, such a configuration can be part of a pair of paths
 $(P',Q')$  in the image of $\Phi$ 
exactly when $a+d\leq i$ and $b+c\leq i$, where $i-1$ is the sum of the two
 coordinates of
$X$. Let us keep  $b$ fixed, and see what that means for $a$ and $c$.
The value of $a$ can be $b,b+1,\ldots, i-1$, so  $a$ can take $i-b$
different values,
 whereas the value of $c$ can be $1,2,\ldots, i-b$ which is
again $i-b$ different possibilities. 
Note in particular that the second set of
values can be obtained from the first by simply subtracting each value 
from $i$. Then the set of all labeled paths
from $(0,0)$ to $A$ is identical to that of paths from $(1,-1)$ to $C$.
In particular, the distributions of the labels of the edges ending in
   $A$, respectively $C$, are identical, even if we also require that they 
end in a horizontal, or in a vertical edge. 

\begin{figure}[ht]
 \begin{center}
  \epsfig{file=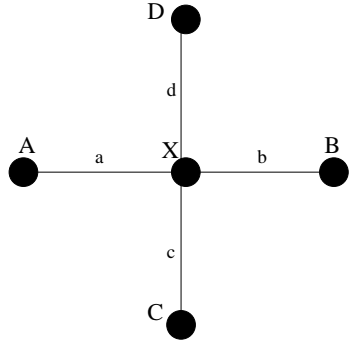}
  \caption{Labels around the point $X$}
  \label{cross}
 \end{center}
\end{figure}

Let $H(X)$ be the set of all pairs of labeled paths 
$((0,0),X) \times ((1,-1),X)$. Now it is easy to see that if any 
labeled path $G$  from $(0,0)$ to $A$ allows $a$ to be in the interval 
$b,b+1,\ldots, i-1$, then the path from $(1,-1)$ to $C$  identical to
$G$ allows $c$
to be in the interval $1,2,\ldots, i-b$. Indeed, the edge preceding $AX$
is either horizontal, and then it must have a label between $b$ and $i-2$,
or it is vertical, and then it must be between 1 and $i-b$ to make it
possible for $a$ to be in the interval 
$b,b+1,\ldots, i-1$. Similarly, if the edge preceding $CX$ is  horizontal,
and it has a label between $b$ and $i-2$, or if it is vertical, and has a label
between 1 and $i-b$, then it makes  it possible for $c$ to be in the interval
$1,2,\ldots, i-b$. (And certainly, if the edge preceding $CX$ is  horizontal,
and it has a label smaller than~$b$, that is good, too). As  the distributions
 of the labels of the edges ending in $A$, respectively $C$ are identical,
 this 
implies that for any fixed values of $b$, there are at least as many
 pairs of paths
in $H(X)$ so that $b+c\leq i$ as there are pairs of paths in $H(X)$ with 
$a\geq b$. (Recall that $i-1$ is the sum of the coordinates of $X$).  In other
words, if the pair $(\alpha,\beta)\in ((0,0),A) \times ((1,-1),C) $ allows
$a\geq b$, then the pair $(\beta,\alpha)\in ((0,0),A) \times ((1,-1),C) $
allows $b+c \leq i$, so we can flip $\alpha$ and $\beta$. We point out that
 this is intuitively not surprising: $a$ has to be {\em at least} a certain
value, while $c$ has to be {\em at most} a certain value, and it is clear
that this second requirement is easier in our labeling. 

By symmetry, if we fix $d$ instead of $b$, the same holds: the number
of  pairs of paths
in $H(X)$  so that $a+d\leq i$ is at least as large as that of 
pairs of paths in $H(X)$ with 
$c\geq d$, and that can be seen again by flipping $\alpha$ and $\beta$.

Finally, this same argument certainly applies  if we want both conditions to be
satisfied: if the pair $(\alpha,\beta)\in ((0,0),A) \times ((1,-1),C)$ allows
$a\geq b$ and $c\geq d$, then the pair $(\beta,\alpha)\in ((0,0),A) \times 
((1,-1),C)$ allows $b+c \leq i$ and $a+d\leq i$. And this is 
what we wanted to prove: there are 
at least as many pairs of paths in $\Pc(n,k)\times \Pc(n,k)$ which are not 
images of $\Phi$ as there are pairs of paths in $\Pc(n,k-1)\times \Pc(n,k+1)$
for which $\Phi$ is not defined. As $\Phi$ is one-to-one, this proves
that $A(n,k-1)\cdot A(n,k+1)\leq A(n,k)^2$, so the sequence $\{A(n,k)\}_k$ is
log-concave for all $n$.

To prove that the sequence $\{U(n,k)\}$ is log-concave, recall that 
half-ascending
permutations in ${\cal U}(n,k)$
  correspond to elements of $\Q(n,k)$, that is, elements of
 $\Pc(n,k)$ in which all edges $a_i$ are horizontal if $i$ is even.
We point out that this implies $B=D$. 
 Then note that 
$\Phi$ does not change the indices of the edges, in other words, if 
$\Phi(P,Q)=(P',Q')$, and a given
edge northeast from $X$  was the $i$th edge of path $P$, then it will be
the $i$th edge of path $Q'$. Therefore, $\Phi$ preserves the property that
all even-indexed edges are horizontal,  so the restriction of $\Phi$
into $\Q(n,k-1)\times \Q(n,k+1)$ maps into $\Q(n,k) \times \Q(n,k)$.
Finally, we need to show that there are more pairs of paths in 
$\Q(n,k) \times \Q(n,k)$
which are not images of $\Phi$ than there are pairs of paths in
$\Q(n,k-1)\times \Q(n,k+1)$ for which $\Phi$ is not defined. Note that
the corresponding fact in the general case was a direct consequence of
the fact that for any labeled path  $((0,0),A)$ was identical to a unique
labeled path   $((1,-1),C)$, and therefore the distributions of the labels
$a$ and $c$ were identical. This remains certainly true if we restrict
 ourselves to paths in which all edges with even indices are horizontal. 
As any restriction of $\Phi$ is certainly one-to-one, this proves
that $U(n,k-1)\cdot U(n,k+1)\leq U(n,k)^2$.
\end{proof}

Now we are in a position to prove the main result of this paper.
\begin{theorem} \label{main} The polynomial $R_n(x)$ has log-concave 
coefficients, for all positive integers $n$.
\end{theorem}
\begin{proof} First suppose that $n$ is even. 
For $n\leq 3$, the statement is true. If $n\geq 4$, then
Lemma \ref{difficult} shows that $R_n(x)=(x+1)^mT_n(x)$. The coefficients
 of $(x+1)^m$
are just the binomial coefficients, which are certainly log-concave
 \cite{sagan}, while the coefficients of $T_n(x)$ are the $U(n,k)$, which
are log-concave by Theorem~\ref{logcon} and the remark thereafter.
 As the product of two polynomials with
log-concave coefficients has log-concave coefficients \cite{stanley}, 
the proof is complete for $n$ even.

If $n$ is odd, then the equivalent of Proposition 
\ref{trivi} is a bit more cumbersome.
Again, we let us make use of symmetry by taking complements, but instead of
 assuming $p_1<p_2$, let us assume that $p_2<p_3$. Taking  $R_{n,m}(x)$ then
 adds the restrictions
 $p_4<p_5$, $p_6<p_7$, $\ldots$, $p_{n-1}<p_n$. Then it is straightforward from
the definition of $t_m(p)$ that $t_m(p)=d(p)$ where $d(p)$ is the number of
 descents of $p$, and we say, for shortness, that the singleton $p_1$ has 0
runs.

So for odd $n$ we have 
$T_n^{odd}(x)=2\cdot\sum_{{p\in S_n}\atop {p_2<p_3}}x^{t_m(p)}=
2\cdot\sum_{{p\in S_n}\atop {p_2<p_3}}x^{d(p)}$, and then, in order to
 see that the coefficients of $T_n^{odd}(x)$ are log-concave, we can repeat
 the argument of Theorem~\ref{logcon}. Indeed, the coefficient of $x^k$ in
 $T_n^{odd}(x)$
equals the cardinality of $\Q'(n,k)$, the subset of $\Pc(n,k)$ in which the
 edges $a_3,a_5,\ldots, a_7$ are horizontal. And the fact that the $|\Q'(n,k)|$
are log-concave can be proved exactly as the corresponding statement
for the $|\Q'(n,k)|=U(n,k)$, that is, by taking the relevant restriction
of~$\Phi$. 

This completes the proof of the theorem for all $n$.
\end{proof}

\newcommand{\journal}[6]{{\sc #1,} #2, {\it #3} {\bf #4} (#5), #6.}
\newcommand{\book}[4]{{\sc #1,} ``#2,'' #3, #4.}
\newcommand{\bookf}[5]{{\sc #1,} ``#2,'' #3, #4, #5.}
\newcommand{\thesis}[4]{{\sc #1,} ``#2,'' Doctoral dissertation, #3, #4.}
\newcommand{\springer}[4]{{\sc #1,} ``#2,'' Lecture Notes in Math.,
                          Vol.\ #3, Springer-Verlag, Berlin, #4.}
\newcommand{\preprint}[3]{{\sc #1,} #2, preprint #3.}
\newcommand{\preparation}[2]{{\sc #1,} #2, in preparation.}
\newcommand{\appear}[3]{{\sc #1,} #2, to appear in {\it #3}}
\newcommand{\submitted}[4]{{\sc #1,} #2, submitted to {\it #3}, #4.}
\newcommand{\JCTA}{J.\ Combin.\ Theory Ser.\ A}
\newcommand{\AdvancesinMathematics}{Adv.\ Math.}
\newcommand{\JournalofAlgebraicCombinatorics}{J.\ Algebraic Combin.}

\newcommand{\communication}[2]{{\sc #1,} #2,personal communication.}


\begin{thebibliography}{99}
\bibitem{wilf}
\communication{H.\ Wilf}
              {Real Zeroes of polynomials that count runs and
               descending runs}
\end{thebibliography}
\end{document}